\documentclass[oneside,a4paper,11pt]{article}
\usepackage{amsfonts, amsbsy, amsmath, amsthm, amssymb, latexsym, verbatim, enumerate,color}
\usepackage{authblk}
\usepackage{graphicx,amssymb,amsmath}
\usepackage{amsthm}
\usepackage{amssymb}
\usepackage{latexsym}
\usepackage{longtable}
\usepackage{epsfig}
\usepackage{hhline}
\usepackage{hyperref}
\usepackage{aliascnt}
\newtheorem{theorem}{Theorem}[section]
\newaliascnt{lemma}{theorem}
\newtheorem{lemma}[lemma]{Lemma} 
\aliascntresetthe{lemma}
\newaliascnt{proposition}{theorem}
\newtheorem{Proposition}[proposition]{Proposition}   
\aliascntresetthe{proposition}
\newaliascnt{corollary}{theorem}
\newtheorem{corollary}[corollary]{Corollary}
\aliascntresetthe{corollary}
\newaliascnt{question}{theorem}

\aliascntresetthe{question}
\newaliascnt{conjecture}{theorem}

\aliascntresetthe{conjecture}

\newtheorem*{theorem*}{Theorem}
\numberwithin{theorem}{section}
\newcommand{\pf}{{\it Proof:\quad}}

\newcommand{\dne}{\hfill $\Box$ \vspace{0.3cm}}

\definecolor{re}{rgb}{1,0.2,0.2}           
 \definecolor{gr}{rgb}{0,1,0}
 \definecolor{bl}{rgb}{0,0,0.6}
 \definecolor{bl2}{rgb}{0,1,0}

\newcommand{\re}{{\mathbb R}}

\newcommand{\I}{{\rm I}}
\newcommand{\el}{{\ell}}

\title{Singular Graphs with Dihedral Group Action}
\author{Ali Sltan Ali AL-Tarimshawy$^1$ \and J. Siemons$^2$}
\date{%
    $^1$Department of Mathematices and Computer ِِApplications, College of Science, Al-Muthanna University, Iraq. \\ alisltan81@yahoo.com \\%
    $^2$School of Mathematics, University of East Anglia, Norwich, Norfolk, NR4 7TJ, UK\\j.siemons@uea.ac.uk \\[2ex]%
       Version of 27 April 2020, compiled \today
}
\begin{document}
\maketitle
{\sc Abstract:}\,   Let $\Gamma$ be a simple undirected graph on a finite vertex set and let $A$ be its adjacency matrix. Then $\Gamma$ is  {\it singular} if  $A$ is singular. The problem of characterising singular graphs is easy to state but very  difficult to resolve in any generality.  In this paper  we investigate the singularity of graphs for which the dihedral group  acts transitively on vertices as a group of automorphisms. 
~\!\!\footnote{\phantom{xxix}{\sc Mathematics Subject Classification:} 05C25, 05E18, 15A18}

\section{\sc Introduction}

Let $\Gamma=(V,E)$ be a graph with vertex set $V$ and edge  set $E.$  Two distinct vertices $v$ and $w$ are  \textit{adjacent}, denoted $v \sim w,$  if $\{v,w\}\in E.$ The \textit{adjacency matrix} of $\Gamma$ is the  matrix whose rows and columns are indexed by the vertices of $\Gamma$ so that the $(v,w)$-entry is equal to $1$ if  $v \sim w,$ and $0$ otherwise. This matrix  is denoted by $A.$  The \textit{characteristic polynomial} of $\Gamma$ is the polynomial $${\rm char}(A,x)={\rm det}(x I_n-A)$$ where $n=|V|$ and  $I_n$ is the corresponding identity matrix. The roots of the characteristic polynomial of $\Gamma$ are the \textit{eigenvalues} of $\Gamma.$ We call  $\Gamma$  \textit{singular} if $A$ is singular. The \textit{spectrum }of $\Gamma$  consists of all eigenvalues  of $A$ and so $\Gamma$ is singular if and only if $0$ belongs to the spectrum of $\Gamma.$  All graphs in this paper are  undirected, simple and finite.

The problem of graph singularity first arose in structural chemistry in the context of H\"uckel Theory~\cite{graovac1972graph}.  In Lagrangian mechanics, when a discrete system is represented by a graph, the eigenvalues of that graph correspond, in some approximations, to the energy levels occurring in the system, further references can be found in~\cite{sieza}. It has been observed that in molecular dynamics the eigenvalue $0$ indicates an instability of the system.  The first mathematical paper on the subject appears to be  Collatz and  Sinogowitz~\cite{von1957spektren} in 1957 who asked for a classification of all finite non-singular graphs. Recently graph singularity has become relevant in other areas of mathematics as well. This applies to Cayley graphs and the representation theory of finite groups in general, see also \cite{liu2018eigenvalues, sieza, Ali}.  There are also  applications in combinatorics and algebraic geometry. For instance, in \cite{muller2005some} M\"{u}ller and  Neunh\"{o}ffer compute the rank of a certain (extremely large) adjacency matrix. The question of whether  this matrix is singular or not has significant importance  for the Foulkes conjecture on the plethysms of symmetric groups. It has become evident that the problem of  Collatz and Singowitz can not be solved in its original version, it would simply not be feasible to determine all singular graphs. 
 
 In this paper we focus on graphs whose automorphism group acts transitively on  vertices. Many results on eigenvalues of graphs of this kind are available, an excellent survey is given by  Liu and Zhou~\cite{liu2018eigenvalues}. In particular, Babai~\cite{Babai} has determined the eigenvalues of graphs with dihedral group action  using the character theory of this group. 
 

The point of view in this article is ring theoretical.  We show that if the graph $\Gamma$ admits a dihedral group  action then certain polynomials $\Psi_{\Gamma}(x)$ can be associated to $\Gamma$ in a natural way. These polynomials then allow us to determine the singularity of $\Gamma$  through divisibility criteria for $\Psi_{\Gamma}(x).$ It is easy to prove that a graph with vertex transitive dihedral action is a Cayley graph. The results therefore apply in particular to Cayley graphs. 
Since any non-solvable group contains dihedral subgroups, the graphs studied here occur as induced subgraphs in the Cayley graphs of any non-sovable group.  Further results on singular graphs can be found in the PhD thesis~\cite{Ali} of the first author.

 To state our results we require the following definitions. Let $\Gamma=(V,E)$ be a graph. If $g$ is  permutation of $V$ and $v\in V$ we write $v^{g}$ for the image of $v$ under $g.$ Then $g$ is an \textit{automorphism} of $\Gamma$ provided that $v\sim u$ if and only $v^{g}\sim u^{g}$ for all $v,\,u\in V.$  The set of all automorphisms of $\Gamma$ is the \textit{automorphism group} of $\Gamma,$ denoted by ${\rm Aut}(\Gamma).$ 
  
 Let $G$ be a finite group. A subset $H$ of $G$ is called a \textit{connecting set} if  (i) $H^{-1}=\{h^{-1}: h \in H\}=H,$  (ii) $1_G \notin H$ and  (iii) $H$ generates $G.$ In this case we define a graph with vertex set $V=G$ by declaring  two vertices $u,v \in G$ to be adjacent  if and only if $vu^{-1} \in H.$ The resulting graph is  the \textit{Cayley graph} of $G$ with connecting set $H,$  denoted by ${\rm Cay}(G,H).$ It is easy to verify that for each $g\in G$ the map $g\!:\,v\to vg$ is an automorphism of $\Gamma$ and further that $G\subseteq {\rm Aut}(\Gamma).$

The dihedral group of order $2n$ is denoted by $D_{n}:=\langle a,b\vert a^n=b^2=1,bab=a^{-1}\rangle;$ its rotation subgroup is denoted by $C_{n}:=\langle a\rangle.$ 
 
 \bigskip
\begin{Proposition}\label{AA1} Let $\Gamma$ be a graph with vertex set $V$ and suppose that the dihedral group $D_{n}$ acts transitively on $V$ as a group of automorphisms of $\Gamma.$  Then either \\[5pt]
(i)\,\, $\Gamma={\rm Cay}(C_{n},H)$ for some connecting set $H\subseteq C_{n},$ or \\[5pt]
(ii)\,\, $\Gamma={\rm Cay}(D_{n},H)$ for some connecting set $H\subseteq D_{n}.$ 
\end{Proposition}

Correspondingly we say that $\Gamma$ is {\it cyclic} of order $n$ if $\Gamma$ is isomorphic to ${\rm Cay}(C_{n},H)$ for some $n,$ and that   $\Gamma$ is {\it dihedral} of order $2n$ if  $\Gamma$ is isomorphic to ${\rm Cay}(D_{n},H)$ for some $n.$ 
Let $\Phi_{d}(x)$ be the $d^{th}$ cyclotomic polynomial. We have 

\bigskip
\begin{theorem}\label{AA2}Let $\Gamma$ be a cyclic graph with of order $n$ with vertex set $V=C_{n}=\{a,a^{2},\dots,a^{n}\}.$ Let $$\Psi_{\Gamma}(x)=m_{2}x+\dots+m_{n}x^{n-1}$$ where $m_{i}=1$ if $a\sim a^{i}$ and $m_{i}=0$ otherwise. Then $\Gamma$ is singular if and only if $\Phi_{d}(x)$ divides $\Psi_{\Gamma}(x)$ for some divisor $d$ of $n.$ 
\end{theorem}

Evidently the coefficients of $\Psi_{\Gamma}(x)$ are determined  by the adjacency matrix of $\Gamma,$ here its first row, or equivalently by the connecting set of the corresponding Cayley graph. The theorem therefore can also be used to define connecting sets $H$ for which ${\rm Cay}(C_{n},H)$ is singular, or non-singular, as required. In general the connection between polynomials and adjacency matrices will be more complicated and this relationship is oulined in Section 2.  

\bigskip
\begin{corollary}\label{AA21} Let $\Gamma$ be cyclic graph with at least one edge. Suppose that the order of $\Gamma$ is prime. Then $\Gamma$ is not singular. 
\end{corollary}

Both results above are also available in a different setting in ~\cite{lal2011non} where additional  properties of cyclic graphs can be found.  
The simplest singular cyclic graph is the $4$-cycle. Indeed, an $n$-cycle is singular if and only if $n$ is divisible by $4,$ see \cite{godsil2013algebraic}. For $n=4$ there are two cyclic singular graphs and there are four singular cyclic graphs if $n=6$ or $n=9.$ The number of singular graphs for $n\leq 11,$ with arbitrary automorphism group, is the sequence  A133206 in the Online Encyclopedia of Integer Sequences~\cite{oeis}.


Next we come to the second class of graphs on which the dihedral group acts vertex transitively. Here  $\Gamma={\rm Cay}(D_{n},H)$ for some connecting set $H$ in $D_{n}.$  We let $V=V'\cup V''$ where $V':=\{a,a^{2},\dots,a^{n}\}=C_{n}$ and $V'':=D_{n}\smallsetminus C_{n}.$ 

\bigskip
\begin{theorem}\label{AA3}Let $\Gamma$ be a dihedral graph of order $2n$ with vertex set $V=V'\cup V''.$ Let $$\Psi'_{\Gamma}(x)=u_{1}+u_{2}x+\dots+u_{n}x^{n-1}$$ where $u_{i}$ is the number of walks of length $2$ from $a$ to $a^{i}$ with intermediate vertex belonging to  $V'.$ Furthermore,  let $$\Psi''_{\Gamma}(x)=w_{1}+w_{2}x+\dots+w_{n}x^{n-1}$$ where $w_{i}$ is the number of walks of length $2$ from $a$ to $a^{i}$ with intermediate vertex belonging to  $V''.$ Then $\Gamma$ is singular if and only if $\Phi_{d}(x)$ divides $\Psi'_{\Gamma}(x)-\Psi''_{\Gamma}(x)$ for some divisor $d$ of $n.$ 
\end{theorem}

We note two corollaries: 

\bigskip
\begin{corollary}\label{pap212}
Let $\Gamma$ be a dihedral graph of order $2n$ with vertex set $V=V'\cup V''.$ Let $v\in V'$ and let $T$ be the set of neighbours of $v.$ Suppose that  $|T\cap V'|=|T\cap V''|.$
Then $\Gamma$ is singular. 
\end{corollary}

\bigskip
\begin{corollary} \label{pap210} With the same notation as in \autoref{pap212} suppose that $n$ is prime and  $|T\cap V'|\neq |T\cap V''|.$ Let $u_{i}$ and $w_{i}$  be as in Theorem~\ref{AA3}. Then $\Gamma$ is singular if and only if there exists $c\in \{\pm 1\}$ so that $u_{i}-w_{i}=c$ for all $i=1,\dots,n$.
\end{corollary}


The proof of these results will be given in Section 4.

\section{\sc Circulant Matrices}

We begin with some results about circulant matrices which will be needed in the paper, most of these can be found in  \cite{daviscirculant}.

A matrix $Y$ of the shape {$$Y = 
 \begin{pmatrix}
  y_1 & y_2 & y_{3}&\cdots & y_{n} \\
  y_{n} & y_1 & y_{2}&\cdots & y_{n-1} \\
  y_{n-1} & y_n & y_{1}&\cdots & y_{n-2} \\
  \vdots  & \vdots  & \vdots & \ddots & \vdots  \\
 y_2 & y_3 & y_{4}&\cdots & y_1
 \end{pmatrix}\,,$$ }\\[5pt] where the elements of each row are those of the previous row, moved one position to the right and wrapped around, is a {\it circulant} matrix.  It is denoted as $Y={\rm circ}(y)$ where $y=(y_1,y_2,...,y_{n})$ is the first row of $Y.$  Thus $Y$ is circulant if and only if $Y_{i,j}=y_{j-i+1}$ when the indices $i,\,j $ are taken modulo $n.$ The set $\mathcal{C}$ of all circulant matrices is a vector space over $\re$ and ${\rm circ}:\re^{n}\to \mathcal{C}$ is a vector space isomorphism. 
 
 In addition, $\mathcal{C}$ is also a commutative ring. To see this denote $W\!:={\rm circ}(0,1,0,\dots,0)$ and let  $\sigma:\re[x]\to \mathcal{C}$ be the substitution map, given by $\sigma(f(x))=f(W).$ This is a surjective ring homomorphism, and since $W^{n}=\I,$ its kernel is the ideal $(x^{n}-1)\re[x].$ Hence \begin{equation}\label{CC1}\re[x]/(x^{n}-1)\,\simeq \,\mathcal{C}\end{equation} as rings. 
 Let $B=(b_{i,j})$ be an arbitrary $n\times n$ matrix. Then the polynomial  \begin{equation}\label{CC2}\Psi(B,x)=b_{1,1}+b_{1,2} x+b_{1,3}x^{2}+\dots+b_{1,n}x^{n-1}  \end{equation} is the {\it associated polynomial} of $B.$ Evidently, $B$ is circulant if and only if $\sigma(\Psi(B,x))=B.$ In particular, we have the following  
 
 \medskip
 \begin{lemma}\label{CC3} Let $B$ and $C$ be $n\times n$ matrices. Then $\Psi(B+C,x)=\Psi(B,x)+\Psi(C,x).$ Furthermore,  if $B$ and $C$ are circulant, then $\Psi(B\cdot C,x)\equiv\Psi(B,x)\cdot\Psi(C,x)\mod (x^{n}-1).$
 \end{lemma}


We will use the fact that the eigenvalues and singularity of a circulant matrix can be determined from its associated polynomial. 
 
Let $\Phi_{\el}(x)$ denote the ${\el}^{th}$ \textit{cyclotomic polynomial}. Thus $\Phi_{\el}(x)$ is the unique irreducible integer polynomial with leading  coefficient $1$ so that $\Phi_{\el}(x)$ divides  $x^{\el}-1$ but  does not divide  of $x^k-1$ for any $k<{\el}$. Its roots are all primitive ${\el}^{th}$  roots of unity. So  $$\Phi_{\el}(x)=\prod_{1\leq m<{\el}}(x-e^\frac{2\pi im}{{\el}})$$ where $gcd(m,{\el})=1.$ 
 
\bigskip

 \begin{lemma}\label{pap24}\cite{kra2012circulant}
Let $Y$ 
 be an $n\times n$ circulant matrix with associated polynomial $\Psi(x)=\Psi(Y,x)$ and let $\Phi_d(x)$ be the $d^{th} $ cyclotomic polynomial. Then $Y$ is singular if and only if $\Phi_d(x)$ divides $\Psi(x)$ for some divisor $d$ of $n.$ 
 Furthermore,
 the nullity of $Y$ is $\sum\varphi(d)$ where the sum is over all $d$ such that $\Phi_d(x)$ divides $\Psi(x)$   and where $\varphi(d)$ is the  Euler totient function.
 \end{lemma}

In some instances these results can be used for other matrices as well. A case in point are matrices  of the form $$X = 
 \begin{pmatrix}
  x_1    & x_2 & x_{3}&\cdots & x_{n} \\
  x_{2} & x_3 & x_{4}&\cdots & x_{1} \\
  x_{3} & x_4 & x_{5}&\cdots & x_{2} \\
  \vdots  & \vdots  & \vdots  &\ddots & \vdots  \\
 x_{n} & x_1 & x_{2}&\cdots & x_{n-1}
 \end{pmatrix}$$ \\[5pt] where the elements of each row are those of the previous row, moved one position to the left and wrapped around. Such a matrix is  called  {\it anti-circulant.}  Clearly, $X$ is anti-circulant if and only if $X_{i,j}=x_{j+i-1}$ when indices are taken modulo $n$ and when $x_{1},...,x_{n}$ is the first row of $X.$  Note, an anti-circulant matrix is automatically symmetric. It is easy to see that the product of two anti-circulant matrices is circulant, see also \cite[Theorem~ 5.1.2]{daviscirculant}. Therefore the singularity of an anti-circulant matrix can be decided from the associated polynomial of its square.

\section{\sc Actions of the Dihedral Group }
We consider the faithful permutation actions of the dihedral group $D_{n}$ of order $2n$ with rotation subgroup $C_{n}$ of order $n.$ We assume throughout $n\geq 3.$

\bigskip
\begin{lemma}\label{pap211} Suppose that  $D_{n}$ acts transitively and faithfully on some set  $V.$ Then  \\[5pt](i) $|V|=n$ and $C_n$ acts regularly on $V,$ or  \\[5pt] (ii)  $|V|=2n$ and $D_{n}$ acts regularly on $V$ while $C_n$ has two orbits on $V.$ 
\end{lemma}

\pf Select some $v\in V$ and let $F\subset G$ be the stabilizer of  $v.$ Put $X:= F\cap C_n.$ Since $C_n$ is cyclic $X$ is the only subgroup of order $|X|$ in $C_{n}.$ In particular, $X$ is normal in $G.$ By a general lemma on permutation groups, for any $g\in G$ we have that $F^g=g^{-1} F g$ is the stabilizer of $v^{g}$ for any $g\in G,$ and so $X=X^g\subseteq F^g$ fixes $v^g$ for all $g.$ By transitivity, $X$ fixes all points of $V$ and  that means $X=1_G,$ as $F$ is faithful on $V.$ In particular, $|F|=1$ or $|F|=2.$  Since $|G|=|V||F|$ by the orbit stabilizer theorem we have the two options in the lemma. \dne

The following proves \autoref{AA1}:

\bigskip
\begin{lemma}\label{pap212c} Let $\Gamma$ be a graph on the vertex set $V$ and suppose that  $D_{n}$ acts transitively on $V$ as a group of automorphisms. Then \\[5pt](i)\,\, $\Gamma={\rm Cay}(C_{n},H)$ where $H$ is a suitable connecting set in $C_{n},$ or \\[5pt](ii)\, $\Gamma={\rm Cay}(D_{n},H)$ where $H$ is a suitable connecting set in $D_{n}.$ \end{lemma}

\pf  By  \autoref{pap211} the graph has a group of automorphisms that acts regularly on its vertices, $C_{n}$ in the first case, and $D_{n}$ in the second case. The result follows from Subidussi's  Theorem, see \cite[Lemma~3.7.2]{godsil2013algebraic}.\dne

First we consider the case when $\Gamma={\rm Cay}(C_{n},H)$ is cyclic. We enumerate its vertices $V=C_{n}$ as $${\rm (E1)}\!:\,\,a, a^{2},\dots, a^{n-1},a^{n}.$$ Then we have 

\begin{lemma}\label{pap212a} Let $A$ be the adjacency matrix of $\Gamma={\rm Cay}(C_{n},H)$ with respect to the enumeration {\rm (E1)} of the vertices of   $\Gamma.$ Then $A$ is circulant. \end{lemma}

\pf To decide the $(i,j)$-entry of $A$ consider the pair $a^{i},\,a^{j}.$ By definition we have $a^{i}\sim a^{j}$ if and only if  $a^{j-i}\in H$ if and only if  $a^{i+k}\sim a^{j+k},$ for all $k.$ Therefore  $A$ is circulant. 
\dne

{\it Proof of \autoref{AA2}:}\, By Lemma~\ref{pap212a} $A$ is circulant with associates polynomial $\Psi_{\Gamma}(x)=\Psi(A,x).$ The result therefore follows from \autoref{pap24}.\dne 

{\it Proof of \autoref{AA21}:}\, Here we have $\Psi_{\Gamma}(x)=m_{2}x+\dots+m_{n}x^{n-1}=x(m_{2}+\dots+m_{n}x^{n-2})$ with $n$ prime. Neither $\Phi_{1}(x)$ nor $\Phi_{n}(x)$ divides  $(m_{2}+\dots+m_{n}x^{n-2})$ and so the result follows from  \autoref{AA2}.\dne 



From now on we consider the case when $\Gamma$ is a dihedral graph. Hence by \autoref{AA1}  we assume that  $\Gamma={\rm Cay}(D_{n},H)$ where $H$ is a connecting set in $D_{n}.$   We enumerate $V=D_{n}$ as  
$${\rm (E2)}\!:\, a, a^{2},\dots, a^{n-1},a^{n};\, ab, a^{2}b,\dots, a^{n-1}b,a^{n}b.$$ We have 

\bigskip
\begin{lemma}\label{pap213a} Let $\Gamma$ be dihedral and let  $A$ be its adjacency matrix, with rows and columns arranged according  to {\rm (E2).}  Then $$A=\left(\begin{array}{cc}M & N \\N& M\end{array}\right)$$ is a blocked matrix with $M$ circulant  and symmetric, of size $n\times n.$ In addition,   $N$ is anti-circulant and symmetric in particular.  Furthermore $MN=NM.$  
\end{lemma}

\pf As $\Gamma$ is an undirected graph  its adjacency matrix is blocked as  $$A=\left(\begin{array}{cc}M_{1}& N_{1} \\N_{2} & M_{2}\end{array}\right)$$ where $M_{1}^{T}=M_{1},$ $M_{2}^{T}=M_{2}$ and $N_{1}^{T}=N_{2}.$ First we show that $M_{1}=M_{2}.$ For this consider the $(i,j)$-entry of these matrices. Observe that $a^{i}\sim a^{j}$ if and only if  $a^{j-i}\in H,$ while $a^{i}b\sim a^{j}b$ if and only if  $a^{j}bba^{-i}=a^{j-i}\in H$ (for E2). It follows that $M_{1}=M_{2}.$  (Notice, $a^{\ell}b=ba^{-\ell}$ for all $\ell.$) To show that $M_{1}$ is circulant note that $a^{i}\sim a^{j}$ if and only if  $a^{j-i}\in H$ if and only if  $a^{i+k}\sim a^{j+k}$  for  all $0\leq k,$ and so $M_{1}=M=M_{2}$ are circulant. 

Next suppose that the vertices are arranged according to \,(E2). Then the $(i,j)$-entry of $N_{1}$ is determined by  $a^{i}\sim a^{j}b$ if and only if $a^{j}ba^{-i}=a^{j}a^{i}b=a^{i+j}b\in H.$ It follows that $N_{1}$ is anti-circulant, and in particular symmetric. Hence $N_{1}=N_{2}.$

To complete the proof  we consider $MN$ and $NM$ where $N=N_{1}=N_{2}.$ The  $(i,j)$-entry of $MN$ is the inner product of the $i^{th}$-row of $M$ with the $j^{th}$-column of $N.$ We have $M_{i,k}=1$ if and only if  $a^{k}a^{-i}=a^{k-i}\in H.$ Similarly, we have $N_{k,j}=1$ if and only if  $a^{k}a^{j}b=a^{k+j}b\in H.$ It follows that $(MN)_{ij}$ is the number of $k$ such that $a^{k-i}\in H$ and $a^{k+j}b\in H.$

To compute the  $(i,j)$-entry of $NM$ consider the $i^{th}$-row of $N$ and the $j^{th}$-column of $M.$ We have $N_{i,\ell}=1$ if and only if  $a^{i}a^{\ell}b=a^{i+\ell}b\in H.$ Furthermore, we have $M_{\ell,j}=1$ if and only if  $a^{-\ell}a^{j}=a^{j-\ell}\in H.$ It follows that $(NM)_{ij}$ is the number of $\ell$ such that $a^{i+\ell}b\in H$ and $a^{j-\ell}\in H.$  We have $a^{k-i}\in H$  if and only if $a^{i-k}\in H$ and set $\ell:=k-i+j.$ From this the condition $a^{i+\ell}b\in H$ and $a^{j-\ell}\in H$ turns into $a^{k+j}b\in H$ and $a^{k-i}\in H.$ It follows that  $MN=NM$ as required.\dne 

When $D_{n}$ acts transitively but not faithfully then the action has a kernel $K\subseteq C_{n}.$ In this case $D_{n}/K$ is again a dihedral  group (allowing for the degenerate case $K=C_{n}$) which now acts faithfully. Hence the above results also cover actions which are not faithful. 

The next step would be to consider the case when the group has several orbits. Some of the above techniques can be adapted to deal with this situation. 

\section{\sc The Main Results}

In this section we consider dihedral graphs of order $2n.$ According to \autoref{pap213a}  the adjacency matrix of $\Gamma$ is
$$(*)\!:\,A=\left(\begin{array}{cc}M & N \\N & M\end{array}\right)$$ where furthermore $MN=NM$ if vertices are arranged as in  (E2).

\bigskip
\begin{theorem}\label{pap213}
Let $\Gamma$ be a dihedral graph of order $2n$ with  adjacency matrix $A$ blocked as in {\rm(*)}.  Then the following hold: \\[5pt](i) \,\,The characteristic polynomial of $\Gamma $ is given by $${\rm char}(A,x)={\rm char}(M+N,x)\cdot {\rm char}(M-N,x).$$ In particular, $\lambda$ is an eigenvalue of $\Gamma$ if and only if $\lambda$ is an eigenvalue of $M+N$ or of $M-N.$ 
\\[5pt]
(ii) \,If the vertices of $\Gamma$ are arranged according to {\rm (E2)} then $\det(A)=\det(M^{2}-N^{2}).$ In particular, $\Gamma$ is singular if and only if $M^{2}-N^{2}$ is singular.
\end{theorem}

{\sc Remark:} In (ii) it is essential that vertices  are indeed arranged as in (E2), the determinant formula does not hold in general.  

{\it Proof:} (i)\, We compute the characteristic polynomial of $A$ as follows:
\begin{eqnarray}\label{pap2_7}
{\rm char}(A,x)&=&{\rm det}(x I_{2n}-A)\nonumber\\
  &=&{\rm det}
  \left( {\begin{array}{cc}
   x I_n-M& -N \\
   -N &x I_n- M \\
  \end{array} } \right) \nonumber\\
&=& {\rm det}\left( {\begin{array}{cc}
  x I_n- M-N & x I_n- M-N \\
   -N& x I_n- M\\
  \end{array} } \right)\nonumber\\
  &=&{\rm det}\left( {\begin{array}{cc}
   x I_n-(M+N) & 0 \\
   N& x I_n- (M-N) \\
  \end{array} } \right)\nonumber\\
  &=& {\rm det}(x I_n-(M+N))\cdot  {\rm det}(x I_n- (M-N))\end{eqnarray}
where $0$ is an $n\times n$ zero matrix. The remaining statement in (i) is  evident. 

(ii) From (i) we have $\det(A)=\det(M+N)\cdot \det(M-N).$ By \autoref{pap213a} we have $MN=NM$ for (E2), hence $(M+N)\cdot(M-N)=M^{2}-N^{2}$ and thus the result. \dne

By \autoref{pap213a} we have that $M$ is circulant and $N$ is anti-circulant. In general it would be very difficult to evaluate ${\rm char}(M+N,x)$ or $ {\rm char}(M-N,x)$ in any generality. 

However, this can be avoided. The essential point is that as $N$ is anti-circulant it follows that $N^{2}$ is circulant, see our earlier comment in Section 2. In particular, $M^2-N^2$ is circulant.

In the following we denote $V':=C_{n}=\{a,a^{2},\dots,a^{n}\},\,\,\,V'':=D_{n}\smallsetminus C_{n}$ and $H':=H\cap V',\,\,\,H'':=H\cap V''.$

\bigskip
\begin{lemma}\label{AAB} Let $\Gamma={\rm Cay}(D_{n},H)$ with $M$ and $N$ as in \autoref{pap213}[ii]. \\[5pt](i)\,\, For $1\leq i\leq n$ let $$u_{i}:=|\{  (h_{1},h_{2})\,\,\big|\,\, h_{2}h_{1}=a^{i-1}\text{\,\,with\,\,}h_{i}\in H'\}|.$$ Then $u_{i}$ is the number of walks of length $2$ from $a$ to $a^{i}$ with intermediate vertex belonging to $V'.$ Furthermore, $$\Psi'(x):=u_{1}+u_{2}x+u_{3}x^{2}+\dots+u_{n}x^{n-1}$$ is the associated polynomial of $M^{2}.$ 
\\[5pt] (ii)\,\, For $1\leq i\leq n$ let $$w_{i}:=|\{  (h_{1},h_{2})\,\,\big|\,\, h_{2}h_{1}=a^{i-1}\text{\,\,with\,\,}h_{i}\in H''\}|.$$ Then $w_{i}$ is the number of walks of length $2$ from $a$ to $a^{i}$ with intermediate vertex in $V''.$ Furthermore,  $$\Psi''(x):=w_{1}+w_{2}x+w_{3}x^{2}+\dots+w_{n}x^{n-1}$$ is the associated polynomial of $N^{2}.$ 
\end{lemma}

\pf (i) The rows and columns of $M$ are indexed by $a,\,a^{2} ,\dots, a^{n}.$ This implies that the $i^{th}$ entry $u_{i}$ in the first row of $MM^{T}$ is the number of vertices with an edge from $a\in V'$ to some $a^{\ell}=h_{1}a\in V'$ so that there is an edge  from $a^{\ell}$ to $a^{i}=h_{2}a^{\ell}\in V'.$ Hence,  $a^{i}=h_{2}h_{1}a$ with $h_{1},\,h_{2}\in H'$ and so $u_{i}$ is the number of pairs $(h_{1},h_{2})$ with  $a^{i-1}=h_{2}h_{1}.$ Using the definition of the associated polynomial the result follows since $M^{T}=M.$ 

(ii) The rows of $N$ are indexed by $a,\,a^{2} ,\dots,a^{n}$ while its columns are indexed by $ab,\,a^{2}b ,\dots, a^{n}b,$ according to (E2). This implies that the $i^{th}$ entry $w_{i}$ in the first row of $NN^{T}$ is the number of vertices with an edge from $a\in V'$ to some $a^{\ell}b=h_{1}a\in V''$ so that there is an edge  from $a^{\ell}b$ to $a^{i}=h_{2}a^{\ell}b\in V'.$ Hence,  $a^{i}=h_{2}h_{1}a$ with $h_{1},\,h_{2}\in H''$ and so $w_{i}$ is the number of pairs $(h_{1},h_{2})$ with  $a^{i-1}=h_{2}h_{1}.$ The result follows since $N$ is anti-circulant, by \autoref{pap213a}.\dne

{\it Proof of \autoref{AA3}:} Since $M^{2}$ and $N^{2}$ are circulant we have 
$\Psi(M^{2},x)-\Psi(N^{2},x)=\Psi(M^{2}-N^{2},x)$ by \autoref{CC3}. By \autoref{pap213}(ii)\, and \autoref{pap24} it follows that $\Gamma$ is singular if and only if $\Psi'(x)-\Psi''(x)$ is divisible by $\Phi_{d}(x)$ for some divisor $d$ of $n.$ \dne 

{\it Proof of \autoref{pap212}:}  Let $v\in V'$ have neighbour set $T.$ From the definition of the polynomials we have $\Psi'(1)=|T\cap V'|^{2}$ and $\Psi''(1)=|T\cap V''|^{2}.$ Hence 
$\Psi'(1)=\Psi''(1)$ and so  $\Phi_1(x)=x-1$ divides $\Psi'(x)-\Psi''(x).$ Thus $\Gamma$ is singular by \autoref{AA3}.\dne

{\it Proof of \autoref{pap210}:} When $n$ is prime we need to consider only two cyclotomic polynomials, $\Phi_{1}(x)=1-x$ and $\Phi_{n}=1+x+x^{2}+\dots+x^{n-1}.$ We have seen above that $1-x$ divides $\Psi'(x)-\Psi''(x)$ if and only if $|T\cap V'|=|T\cap V'|.$ Therefore it is sufficient to consider the second case. The condition in the Corollary is equivalent to saying that $\Psi'_{\Gamma}-\Psi''_{\Gamma}=c\cdot \Phi_{n}.$ The result now follows from  \autoref{AA3}.\dne

\end{document}